\documentclass[12pt]{amsart}
\usepackage{amssymb,amsmath}
\usepackage{andy}
\usepackage{amsthm}
\usepackage{pdfsync}

\newtheorem{thm}{Theorem}[section]
\newtheorem{prop}[thm]{Proposition}

\theoremstyle{definition}

\theoremstyle{remark}

\renewcommand{\v}{v_{(m)}}

\newcommand{\Sk}{\mathcal S_{_{\!(k)}}}
\newcommand{\Sm}{\mathcal S_{_{\!(2m+1)}}}
\newcommand{\X}{X_{(m)}}
\newcommand{\Y}{Y_{(m)}}

\begin{document}

\title[Higher order exterior derivative] {Higher order analogues of  exterior derivative
\footnote{This is an electronic reprint of the original article published in the Bulletin of the  Institute of Mathematics of the Academia Sinica (New Series) {\bf 8 (3)} (2013) 389 -- 398. This reprint differs from the original in pagination and typographical detail.}}
\author[Lanzani]{Loredana Lanzani$^*$}
\thanks{To appear in Bull. IMAS, New Series.}
\thanks{$^*$ Supported by a National Science Foundation IRD plan, and in part by award DMS-1001304.}

\address{National Science Foundation and University of Arkansas\\ Fayetteville AR 72701}

\address{Dept. of Mathematics \\
University of Arkansas Fayetteville, AR 72701}

\email{loredana.lanzani@gmail.com}
\maketitle

\begin{abstract}
 We give
   new examples of linear differential operators of order $k=2m+1$ (any given odd integer)  that are invariant under the isometries of $\mathbb R^n$ and satisfy so-called $L^1$-duality estimates and
div/curl inequalities. 
\end{abstract}

\section{Introduction}\label{S:introduction}
The purpose of this note is to exhibit (elementary) 
examples of $k$th-order linear differential operators
$\{\Sk\}_{_{\! k}}$ acting 
 on $\mathbb R^n$ that 
 can be regarded as higher order analogues of the exterior derivative 
  complex
 $$
 d: C^{\infty, c}_q(\mathbb R^n)\ \to\ C^{\infty, c}_{q+1}(\mathbb R^n),\qquad 0\leq q\leq n
 $$
(Here $C^{\infty, c}_q(\mathbb R^n)$ and $C^{\infty, c}_{q+1}(\mathbb R^n)$ stand for the
 $q$-forms and $(q+1)$-forms on $\mathbb R^n$ whose coefficients are smooth and compactly supported.)
 More precisely we require that, for each $k$,
 $\Sk$ map $q$-forms to $(q+1)$-forms and $\Sk\circ\Sk =0$; that the Hodge Laplacian for $\Sk$, namely the operator
  $\Sk\Sk^*+\Sk^*\Sk$,  be elliptic, and that the first-order operator in this family be the exterior
  derivative (that is, $\mathcal S_1=d$).
We  also require that $\Sk$ and $\Sk^*$ have non-trivial invariance properties
 and
satisfy so-called $L^1$-duality estimates as well as div-curl inequalities (more on these below).
 While 
 various operators satisfying one or more of these conditions were recently constructed for any 
 order 
 $k=1, 2, 3,\ldots$, see \cite{BB3}, \cite{LR} and \cite{VS3}-\cite{VS5}, 
 those
 operators
 fail to be invariant under pullback by the rotations of $\mathbb R^n$ as soon as $k\geq 2$. By contrast, 
 here we 
 define linear differential operators $\Sk$
    of {\em odd} order  
  $$
 k=2m +1,\quad m=0, 1, 2,\ldots ,
 $$
that have the
 same invariance properties as the codifferential $d^*$ (the $L^2$-adjoint of exterior derivative)
  as soon as $k\geq 3$ (i.e. $m\geq 1$);
 that is
 $$
 \Sk\circ\, \psi^*= \psi^*\circ\Sk\qquad \mbox{and}\qquad
 \Sk^*\circ\, \psi^*= \psi^*\circ\Sk^*
 $$
 for any isometry $\psi:\mathbb R^n\to\mathbb R^n$ (as customary, 
 $\psi^*$ denotes the pullback of $\psi$ acting on $q$-forms).
 While such invariance 
  is non-trivial, it is far weaker than the invariance
 of $d$,  which indeed is what should be expected of any linear differential operator of order greater than 1, see 
 \cite[Note 4]{P} and \cite{Te}.
 \\
 
 Specifically, given $m=0, 1, 2, 3, \ldots$,
 we define
 \begin{equation}\label{D:odd-op}
 \Sm := d\,(d^*d)^m\quad \mbox{and, consequently},\quad 
  \Sm^* = (d^*d)^md^*
 \end{equation}
 \noindent It is clear that $\mathcal S_{(1)}= d$ and, more generally, that
 $\Sm$ takes $q$-forms to $(q+1)$-forms and  $\Sm\circ\,\Sm =0$.
 It is also clear that the Hodge Laplacian for $\Sm$ is
 $$
 \square_{(2m +1)}\ =\ \square^{2m+1}\ =\ \square\circ\square\circ\cdots\circ\square\quad 
 $$
  where the composition above is performed $(2m+1)$-many times and
 $$
 \square = dd^* + d^*d
 $$
 is the Hodge Laplacian for the exterior derivative,
    so in particular $\square_{(2m+1)}$ is
 elliptic because it is the composition of elliptic operators \cite{Wo}. \\
 
 Note, however, that
 $$
 d\circ \Sm=0\quad \mbox{and}\quad  d^*\circ\Sm^*=0
 $$
 see \eqref{D:odd-op}, and so the natural compatibility conditions for the data of the Hodge system for $\Sm$ and $\Sm^*$
 are the same as for the system for $d$ and $d^*$.
 As a consequence, the $L^1$-duality inequalities that are relevant to the Hodge system for 
 $\Sm$ and $\Sm^*$ are the same as in \cite[page 61]{LS} and \cite{VS1}, namely
 \begin{prop}[\cite{LS}]\label{P:LS} There is $C=C(n)$ such that
for any  $0\leq q\leq n-2$ and for any
 $f\in C^{\infty, c}_{q+1}(\mathbb R^n)$
  %Then there is $C=C(n)$ such that
  \begin{equation}\label{E:LS-d}
 df=0\quad \Rightarrow\quad  |\langle f, h\rangle|\,\leq\, C \|f\|_{L^1_{q+1}(\mathbb R^n)}\|\nabla h\|_{L^n_{q+1}(\mathbb R^n)}
  \end{equation}
  for any $h\in L^\infty_{q+1}(\mathbb R^n)$ such that $\nabla h\in L^n_{q+1}(\mathbb R^n)$.\\
  
  There is $C=C(n)$ such that for any
    %Suppose that 
    $2\leq q\leq n$ and for any
 $g\in C^{\infty, c}_{q-1}(\mathbb R^n)$
 %. Then there is 
 
% \noindent $C=C(n)$ such that
 \begin{equation}\label{E:LS-d-adj}
 d^*g=0\quad\Rightarrow\quad |\langle g, h\rangle|\,\leq\, C \|g\|_{L^1_{q-1}(\mathbb R^n)}\|\nabla h\|_{L^n_{q-1}(\mathbb R^n)}
  \end{equation}
  for any $h\in L^\infty_{q-1}(\mathbb R^n)$ such that $\nabla h\in L^n_{q-1}(\mathbb R^n)$.
 \end{prop}
 
 Here $L^p_{q\pm 1}(\mathbb R^n)$ denote the spaces of $(q\pm 1)$-forms whose coefficients are in the Lebesgue class $L^p(\mathbb R^n)$, and $\langle \cdot, \cdot\rangle$ denotes the inner product in $L^2_{q\pm1}(\mathbb R^n)$:
 $$
 \langle f, h\rangle = \int\limits_{\mathbb R^n}\!\!f\wedge *h 
  $$
where $*$ denotes the Hodge-star operator for $\mathbb R^n$.\\
 
 We take this opportunity to point out that these inequalities can be restated in the seemingly more invariant, in fact equivalent, fashion (see also \cite[Theorem 1'']{BB3})
 \begin{prop}\label{P:LL} There is $C= C(n)$ such that for any
% Suppose that 
 $0\leq q\leq n-2$ and for any
 $f\in C^{\infty, c}_{q+1}(\mathbb R^n)$
  %Then there is $C=C(n)$ such that
  \begin{equation}\label{E:LL-d}
 df=0\quad \Rightarrow\quad  |\langle f, h\rangle|\,\leq\, C \|f\|_{L^1_{q+1}(\mathbb R^n)}\|d^*h\|_{L^n_{q}(\mathbb R^n)}
  \end{equation}
  for any $h\in L^\infty_{q+1}(\mathbb R^n)$ such that $d^* h\in L^n_{q}(\mathbb R^n)$.\\
  
    There is $C=C(n)$ such that for any 
    %Suppose that
     $2\leq q\leq n$ and for any
 $g\in C^{\infty, c}_{q-1}(\mathbb R^n)$
 %. Then there is 
 
% \noindent $C=C(n)$ such that
 \begin{equation}\label{E:LL-d-adj}
 d^*g=0\quad\Rightarrow\quad |\langle g, h\rangle|\,\leq\, C \|g\|_{L^1_{q-1}(\mathbb R^n)}\|d h\|_{L^n_{q}(\mathbb R^n)}
  \end{equation}
  for any $h\in L^\infty_{q-1}(\mathbb R^n)$ such that $d h\in L^n_{q}(\mathbb R^n)$.
 \end{prop}

We show below that this result is equivalent to each of
 the following div/curl-type inequalities (one for any choice of $m=0, 1, 2, \ldots$) which are  proved with the methods of 
 \cite{LS}:
\begin{thm}\label{T:div-curl-odd} Fix $0\leq q\leq n$ and let $f\in L^1_{q+1}(\mathbb R^n)$ with $df=0$, and $g\in L^1_{q-1}(\mathbb R^n)$ with $d^*g =0$ be given.
Then, for any $m=0, 1, 2, 3,\ldots $, the (unique) $q$-form $\v$ that solves
the system
\begin{equation}
 \label{E:Hodge-odd}
 \bigg\{
 \begin{array} {lll}
  \Sm \v \ =& f &
 \\ 
  \Sm^*\, \v\  \ =& g &
  \end{array}
 \end{equation}
belongs to the Sobolev space $W^{2m, r}_q(\mathbb R^n)$ with $r=n/(n-1)$ whenever $q$ is neither 1 (unless $g=0$) nor $n-1$ (unless $f=0$), and we have

   \begin{equation}\label{E:div-curl-odd}
 \|\v\|_{W^{2m, r}_q(\mathbb R^n)}\leq C
 \big(\|f\|_{ L^1_{q+1}(\mathbb R^n)} +
 \|g\|_{ L^1_{q-1}(\mathbb R^n)}\big).
  \end{equation}
 \end{thm}
  Here $W^{2m, r}_q(\mathbb R^n)$ denotes the space of $q$-forms whose coefficients belong to the Sobolev space $W^{2m, r}(\mathbb R^n)$ of functions that are $2m$-many times differentiable in the sense of distributions and whose derivatives of any order ($0\leq |\alpha|\leq 2m$) are in the Lebesgue class $L^r(\mathbb R^n)$. 
 \begin{prop}\label{P:div-curl-H-1}
 With same hypotheses as Theorem \ref{T:div-curl-odd}, if $q=1$ and $g\neq 0$ a substitute of \eqref{E:div-curl-odd} holds with $\|g\|_{L^1(\mathbb R^n)}$ replaced by $\|g\|_{H^1(\mathbb R^n)}$, where
 $H^1(\mathbb R^n)$ is
 the real Hardy space. 
 If  $q=n-1$ and $f\neq 0$, then \eqref{E:div-curl-odd} holds with 
 $\|f\|_{H^1_{n}(\mathbb R^n)}$ in place of $\|f\|_{L^1_{n}(\mathbb R^n)}$, where $H^1_n(\mathbb R^n)$ is the space of $n$-forms whose coefficients are in $H^1(\mathbb R^n)$.
 \end{prop}
 In the case when $m=0$, Theorem \ref{T:div-curl-odd} and Proposition \ref{P:div-curl-H-1} were proved in \cite{LS}, as in such case we have $\mathcal S_{(1)}= d$ and $W^{0, r}_q(\mathbb R^n) = L^r_q(\mathbb R^n)$, and so Theorem
 \ref{T:div-curl-odd} and Proposition \ref{P:div-curl-H-1} can be viewed as a generalization (actually, as we will see, a consequence) of those earlier results.\\
 
We remark in closing that one could also consider the operators
$$
\mathcal S_{(2m)} := (dd^*)^m\quad \mbox{and}\quad \widetilde{\mathcal S}_{(2m)}:= (d^*d)^m
$$
but these fail to map $q$-forms to $(q+1)$-forms and do not form a complex 
%(and do not %recapture exterior derivative for a special choice of $m$) 
and as such are of not pertinent to this note. 
\bigskip

{\em Aknowledgements.}  I would like to thank R. Palais and C.-L. Terng for helpful discussions. Part of this work was developed while I was visiting the Institute of
Mathematics at the Academia Sinica: I am very grateful for the support and for the kind hospitality.

\section{Proofs}\label{S:Proofs}
We begin by recalling the elliptic estimates for $\square^s=\square\circ\cdots\circ\square$, see \cite{CZ} and e.g.,
\cite{Wo}, \cite{SR}.
\begin{thm}\label{T:ellipt-ests}
Given any $s\in\mathbb Z^+$, we have that
\begin{equation*}
\square^{s}\!:  C^{\infty, c}_q(\mathbb R^n)\to  C^{\infty, c}_q(\mathbb R^n)
\end{equation*}
is invertible, and
\begin{equation}\label{E:big-box-elliptic estimates}
\|(\square^{s})^{-1}\,u\|_{W_q^{2s, r}(\mathbb R^n)}\lesssim
\|u\|_{L_q^{ r}(\mathbb R^n)}
\end{equation}
for any $1<r<\infty$.
\end{thm}
{\em Proof of Theorem  \ref{T:div-curl-odd}}. 
The case $m=0$ was proved in \cite{LS} and here we will show that the estimates 
in the case when $m\in\mathbb Z^+$ follow from the inequalities for $m=0$. Without loss of generality we may assume: 
$f\in C^{\infty, c}_{q+1}(\mathbb R^n)$ and 
$g\in C^{\infty, c}_{q-1}(\mathbb R^n)$, so that each of 
$d^*f$ and $dg$ has smooth and compactly supported coefficients.

 Applying the codifferential $d^*$ to the first equation in 
 \eqref{E:Hodge-odd} 
 and the exterior derivative
 $d$ to the second equation, and then adding the two equations, see \eqref{D:odd-op},
 we find that
  \begin{equation}\label{E:aux-1}
 \square^{m+1}\v = d^*f + dg
 \end{equation}
 Comparing $\v$ with the solution $u$ of the Hodge system for $d$ and $d^*$ with {\em same} data as 
  \eqref{E:Hodge-odd},
  namely
 \begin{equation} \label{E:Hodge-d}
 \bigg\{
 \begin{array} {lll}
 d\,\! u \ \ \,=& f &
  \\
 d^*\, u\ =& g&
 \end{array}
 \end{equation}
we find
$$
\square^{m}\v =  u
$$
 and so the elliptic estimate \eqref{E:big-box-elliptic estimates} (with $s:=m$) grants
 \begin{equation}\label{E:ellipt-m}
 \|\v\|_{W_q^{2m, r}(\mathbb R^n)}\lesssim
\|u\|_{L_q^{ r}(\mathbb R^n)}
 \end{equation}
 for any $1<r<\infty$. 
 On the other hand, by \cite{LS} we have that $u\in L^r_q(\mathbb R^n)$ with $r:=n/(n-1)$ and
 \begin{equation}\label{E:div-curl-LS}
 \|u\|_{L_q^{ r}(\mathbb R^n)}\leq C(n) 
 \big(\|f\|_{ L^1_{q+1}(\mathbb R^n)} + 
 \|g\|_{ L^1_{q-1}(\mathbb R^n)}\big).
\end{equation}
The desired conclusion \eqref{E:div-curl-odd} now follows by combining \eqref{E:ellipt-m} and
\eqref{E:div-curl-LS}.
\qed\\

{\em Proof of Proposition  \ref{P:div-curl-H-1}}. 
The case $m=0$ was proved in \cite{LS} and here we will again only consider $m\in\mathbb Z^+$. As before, we may assume: 
$f\in C^{\infty, c}_{q+1}(\mathbb R^n)$ and 
$g\in C^{\infty, c}_{q-1}(\mathbb R^n)$.
Now \eqref{E:ellipt-m} holds as before,
 and if $q=1$ and $g\neq 0$
it was proved in \cite{LS} that a substitute of \eqref{E:div-curl-LS} holds with $\|g\|_{L^1(\mathbb R^n)}$ replaced by $\|g\|_{H^1(\mathbb R^n)}$, so the proof of Proposition \ref{P:div-curl-H-1} in the case $q=1$ follows by combining \eqref{E:ellipt-m} and the $H^1$-substitute for \eqref{E:div-curl-LS}.
 (The case $q=n-1$ and $f\neq 0$ is proved in a similar fashion.)\qed \\

Next we show that Theorem \ref{T:div-curl-odd} (for any choice of $m=0, 1, 2, \ldots$) is equivalent to Proposition \ref{P:LL}.
\smallskip

{\em Theorem  \ref{T:div-curl-odd} $\Rightarrow$ Proposition \ref{P:LL}}.  To prove \eqref{E:LL-d}, it again suffices to consider the case when $f$ and $h$ have smooth and compactly supported coefficients; given $f$ as in \eqref{E:LL-d} we consider the solution $\v$ (for $m$ fixed arbitrarily) of the system (6) with $g:=0$, namely

\begin{equation*}
 \bigg\{
 \begin{array} {lll}
 d\,(d^*d)^m\, \v \ =& f &
 \\
 \!(d^*d)^md\, \v\  \ =& 0 &
  \end{array}
 \end{equation*}

see \eqref{D:odd-op}, so that
$$
\langle f, h\rangle = \langle d\,(d^*d)^m\, \v, h\rangle\, 
$$
Integrating by parts the right-hand side of this identity
we obtain
$$
|\langle f, h\rangle|= |\langle \v, (d^*d)^md^*h\rangle|
$$
H\"older inequality for $W^{2m,\, n/(n-1)}_q(\mathbb R^n)$ and its conjugate
space 

\noindent $W^{-2m,\, n}_q(\mathbb R^n)$
now grants
$$
|\langle f, h\rangle|\leq \|\v\|_{W_q^{2m,\, n/\!(n-1)}(\mathbb R^n)}
\|(d^*d)^md^*h\|_{W_q^{-2m,\, n}(\mathbb R^n)}
$$
and by Theorem \ref{T:div-curl-odd} it thus follows that
$$
|\langle f, h\rangle|\leq \|f\|_{L^1_{q+1}(\mathbb R^n)}
\|(d^*d)^md^*h\|_{W_q^{-2m,\, n}(\mathbb R^n)}
$$
On the other hand, we have
$$
\|(d^*d)^md^*h\|_{W_q^{-2m,\, n}(\mathbb R^n)}=
\sup\limits_
{_{\|\zeta\|_{W^{2m,\,n/\!(n-1)}_{q}
}\leq 1}}
|\langle (d^*d)^md^*h, \zeta\rangle|
 $$
 Integrating the latter by parts $2m$-many times and applying H\"older inequality for
  $L^n_q(\mathbb R^n)$ and its dual space $L^{n/\!(n-1)}_q(\mathbb R^n)$ we find
  $$
  |\langle (d^*d)^md^*h, \zeta\rangle|\leq
  \|d^*h\|_{L^n_q}\,\|(d^*d)^m\zeta\|_{L^{n/\!(n-1)}_q}
  $$
  but
  $$
  \|(d^*d)^m\zeta\|_{L^{n/\!(n-1)}_q}\leq \|\zeta\|_{W^{2m,\,n/\!(n-1)}_{q}}
  $$
  which concludes the proof of  \eqref{E:LL-d}. To prove \eqref{E:LL-d-adj} it suffices to apply \eqref{E:LL-d}
  to $f:= *h\in C^{\infty, c}_{\widetilde q +1}(\mathbb R^n)$  with $\widetilde q:= n-q$ (recall that $d^*\approx *d*$ and that $*:L^1_q(\mathbb R^n)\to L^1_{n-q}(\mathbb R^n)$ is an isometry).\qed\\
  
  {\em  Proposition \ref{P:LL} $\Rightarrow$ Theorem  \ref{T:div-curl-odd} for any $m=0, 1, 2, \ldots$}
Without loss of generality we may assume, as before, that
$f\in C^{\infty, c}_{q+1}(\mathbb R^n)$ and
$g\in C^{\infty, c}_{q-1}(\mathbb R^n)$.
Fix $m\in \{0, 1, 2, 3, \ldots\}$ arbitrarily and write
$$
\v=\X + \Y
$$
where
\begin{equation}\label{E:system-F}
\left\{
\begin{array}{clr}
 d(d^*d)^m \X&=&f\\
 (d^*d)^md^*\,\X&=&0
\end{array}
\right.
\end{equation}
and 
\begin{equation}\label{E:system-G}
\left\{
\begin{array}{clr}
 d(d^*d)^m\Y&=&0\\
(d^*d)^m d^*\,\Y&=&g
\end{array}
\right.
\end{equation}
see \eqref{D:odd-op}. We claim that
\begin{equation}\label{E:F-est}
\|\X\|_{ W^{2m,\, n/(n-1)}_q}\leq C \|f\|_{ L^1_{q+1}},
\end{equation}
and 
\begin{equation}\label{E:G-est}
\|\Y\|_{ W^{2m, n/\!(n-1)}_q}\leq C \|g\|_{ L^1_{q-1}}
\end{equation}
Note that if $\Y$ solves \eqref{E:system-G} then $\X:= *\Y$ solves 
\eqref{E:system-F} with $f:= *g\in C^{\infty, c}_{\widetilde q +1}(\mathbb R^n)$  and $\widetilde q:= n-q$, and so it suffices to prove 
\eqref{E:F-est} for $f$ and $\X$ as in \eqref{E:system-F}. (Note that the proof of \eqref{E:F-est} is non-trivial only for $q\neq n$, and the hypotheses of Theorem \ref{T:div-curl-odd} require $q\neq n-1$, so all together we may assume $0\leq q\leq n-2$.)
By duality, proving \eqref{E:F-est} is equivalent to showing
\begin{equation}\label{E:F-est-aux}
\big|\langle D^{\beta} \X, \varphi\rangle\big|\leq
C\|f\|_{ L^1_{q+1}} \|\varphi\|_{ L_q^{n}}
\end{equation}
for any for any $\varphi\in C^{\infty, c}_q(\mathbb R^n)$ and for any multi-index $\beta$ of length $s$ (that is, $\beta=(\beta_1,\ldots,\beta_n)\in \mathbb N^n,\ \beta_1+\cdots +\beta_n=s$)
 and for any $0\leq s\leq 2m$, where we have set
$$
D^\beta \X :=\sum\limits_{|I|=q}
\bigg(\frac{\dee^s X_{{(m)}_I}}{\dee x^\beta}\bigg)
dx^I\, .
$$
To this end, write $\varphi = \square^{m+1}\Phi$ for some $\Phi\in C^{\infty, c}_q(\mathbb R^n)$, see Theorem \ref{T:ellipt-ests}; then
$$
\big|\langle D^{\beta} \X, \varphi\rangle\big| = 
\big|\langle D^{\beta} \X, \square^{m+1}\Phi\rangle \big| 
$$
Integrating the right-hand side of this identity parts we find
$$
\big|\langle D^{\beta} \X, \varphi\rangle\big| = 
\big|\langle \square^{m+1}\X, D^\beta\Phi \rangle\big|
$$
But $\square^{m+1}\X = d^*f$, see \eqref{E:system-F} and so
$$
\big|\langle D^{\beta} \X, \varphi\rangle\big| = 
\big|\langle d^*f, D^\beta\Phi \rangle\big| = \big|\langle f, dD^\beta\Phi\rangle\big|.
$$
Applying Proposition \ref{P:LL} to $h:= dD^\beta\Phi\in C^{\infty, c}_{q+1}(\mathbb R^n)$ we conclude
$$
\big|\langle D^{\beta} \X, \varphi\rangle\big|\leq C(n) \|f\|_{L^1_{q+1}}\|d^*dD^\beta\Phi\|_{L^n_q}
\leq C(n) \|f\|_{L^1_{q+1}}\|\Phi\|_{W^{2(m+1), n}_q}
$$
On the other hand, since  we had chosen $\Phi = (\square^{m+1})^{-1}\varphi$,  Theorem \ref{T:ellipt-ests} grants
$$
\|\Phi\|_{W^{2(m+1), n}_q}\lesssim \|\varphi\|_{L^n_q}
$$
which combines with the previous estimates to give the desired inequality.\qed\\

It should by now be clear that Propositions \ref{P:LS} and \ref{P:LL} are equivalent to one another: on the one hand, it is obvious that Proposition \ref{P:LL} $\Rightarrow$ Proposition
\ref{P:LS} (because $\nabla h\in L^n_{q\pm1}\Rightarrow dh\in L^n_{(q+1)\pm 1}$ and 
$d^*h\in L^n_{(q-1)\pm 1}$ and, moreover, $\|dh\|,\ \|d^*h\|\leq \|\nabla h\|$). On the other hand,
it was proved in \cite[page 61]{LS} that Proposition \ref{P:LS} $\Rightarrow$ Theorem \ref{T:div-curl-odd} in the case $m=0$ which in turn, as we have just seen, gives Theorem \ref{T:div-curl-odd} for arbitrary $m$ as well as Proposition \ref{P:LL}.


\begin{thebibliography}{mybib}

\bibitem[A]{A} Adams R. A. {\em Sobolev Spaces}, Academic Press, New York
(1975).

\bibitem[Am]{Am} Amrouche, C. and Nguyen H.-H. {\em New estimates for the div-curl operators and elliptic problems with $L^1$-data in the whole space and in the half space}
J. Diff. Eq. {\bf 250}, 3150-3195 (2010).

\bibitem[BB1]{BB1} Bourgain J., Brezis H.,
{\em On the equation div$Y=f$ and application 
to control of phases}, J. Amer. Math. Soc. {\bf 16}, 393-426 (2003).

\bibitem[BB2]{BB2} Bourgain J., Brezis H.,
{\em New estimates for the Laplacian, the div-curl and
related Hodge systems}, C. R. Math. Acad. Sci. Paris {\bf 338}, 539-543 (2004).

\bibitem[BB3]{BB3} Bourgain J., Brezis H.,
{\em New estimates for elliptic equations and Hodge
 type systems}, J. Eur. Math. Soc. {\bf 9}, 277-315 (2007).

\bibitem[BF]{BF} Baldi A. and Franchi B., {\em Sharp a-priori estimates for div-curl systems in Heisenberg groups}, preprint (2013).

\bibitem[BV]{BV} Brezis, H. and Van Schaftingen, J. {\em Boundary estimates for elliptic systems
with $L^1$-data}, Cal. Var. PDE {\bf 30} 369 - 388 (2007).

\bibitem[CZ]{CZ} Calder\'on A. and Zygmund A. {\em on the existence of certain singular
integrals} Acta Math. {\bf 88},  85 - 139 (1952).

\bibitem[CV]{CV} Chanillo S. and van Shaftingen J. {\em Subelliptic Bourgain-Brezis estimates on groups} Math. Res. Lett. {\bf 16}, 235 - 263 (2009).

\bibitem[G]{G} Gagliardo E. {\em Ulteriori propriet\`a di alcune classi di funzioni di pi\`u variabili}
Ricerche Mat. {\bf 8},  24 - 51 (1959).

\bibitem[HP-1]{HP-1} Hounie J. and T. Picon {\em Local Gagliardo-Nirenberg estimates for elliptic
systems of vector fields}, Math. Res. Lett. \textbf{18} (2011), no.04, 791--804.

\bibitem[HP-2]{HP-2} Hounie J. and T. Picon {\em Local $L^1$ estimates for elliptic systems of complex vector fields}, Proc. AMS, to appear.

\bibitem[KPV]{KPV} Kenig C., Pipher J. and Verchota G. {\em Area integral estimates for higher order equations and systems}, Ann. Inst. Fourier (Grenoble) {\bf 471425 -- 1461}, (1997).

\bibitem[LR]{LR} Lanzani L. and Raich A. S. 
{\em On Div-Curl for Higher Order}, to appear in {\em Advances in Analysis: The Legacy of Elias M. Stein}, Princeton University Press.

\bibitem[LS]{LS} Lanzani L. and  Stein E. M.,
{\em A note on div-curl inequalities}, Math. Res. Lett. {\bf 12}, 57-61 (2005).

\bibitem[M]{M} Mazya, V. {\em Estimates for differential operators of vector analysis involving $L^1$-norm}, J. Math. Soc. {\bf 12}, 221-240 (2010).
 
 \bibitem[Mi]{Mi} Mironescu P. {\em On some inequalities of Bourgain, Brezis, Maz'ya, and
 Shaposhnikova related to $L^1$ vector fields} C.R. Math. Acad. Sci. Paris {\bf 348}, 513-515 (2010).
 
 \bibitem[MM]{MM} Mitrea M. and Mitrea I. {\em A remark on the regularity of the div-curl system} Proc. Am. Math. Soc. {\bf 137}, 1729 - 1733 (2009).
 
\bibitem[Mu]{Mu} Munkres, J. R., {\em Analysis on manifolds}, Westview Press (1990).{\bf}

\bibitem[N]{N} Nirenberg L. {\em On elliptic partial differential equations} Ann. Scuola Norm. Sup. Pisa,
{\bf 13}, 115 - 162 (1959).

\bibitem[P]{P} Palais R. {\em Natural operators on differential forms} Trans. AMS, {\bf 92} no. 1, 125-141 (1959).

\bibitem[SR]{SR} Saint Raymond, X. {\em Elementary introduction to the theory of 
pseudodifferential operators} CRC Press, Boca Raton (1991).

\bibitem[S]{S} Stein E. M. 
{\em Harmonic Analysis} Princeton U. Press (1993).

\bibitem[T]{T} Taylor M. E. {\em Partial Differential Equations III} Springer Verlag (1996).

\bibitem[Te]{Te} Terng C.-L. {\em Natural vector bundles and natural differential operators} Amer. J. Math., {\bf 100}, no. 4, 775-828 (1978).

\bibitem[VS1]{VS1} Van Schaftingen J., 
{\em Estimates for $L^1$-vector fields}, C. R. Math. Acad. Sci. Paris {\bf 338},
23-26 (2004).

\bibitem[VS2]{VS2} Van Schaftingen J.,
{\em Estimates for $L^1$ vector fields with a second order condition} 
C. R. Math. Acad. Sci. Paris {\bf 339},
181-186 (2004).

\bibitem[VS3]{VS3} Van Schaftingen J.,
{\em Estimates for $L^1$-vector fields under higher order 
differential conditions}, J. Eur. Math. Soc. {\bf 10}, 867-882 (2008).

\bibitem[VS4]{VS4} Van Schaftingen J.,
{\em Limiting Sobolev inequalities for vector fields and canceling linear differential operators}, J. Eur. Math. Soc. {\bf 15}, 877-921 (2013).

\bibitem[VS5]{VS5} Van Schaftingen J. {\em Limiting fractional and Lorentz space estimates of differential forms} Proc. Am. Math. Soc. {\bf 138}, 235 - 240  (2010).

\bibitem[Y]{Y} Yung P.-L. {\em Sobolev inequalities for $(0, q)$-forms on CR manifolds of finite type} Math. Res. Lett. {\bf 17} 177 - 196 (2010)

\bibitem[Wo]{Wo} Wong, M.-W.,{\em An introduction to pseudodifferential operators, 2nd. Ed.}, 
World Scientific Co., Singapore (1999)

\end{thebibliography}
\end{document}